\documentclass[a4paper]{amsart}
\usepackage[english]{babel}
\usepackage[utf8]{inputenc}
\usepackage[T1]{fontenc}
\usepackage{amsfonts,amsmath,amssymb,amsthm}
\usepackage{mathtools}
\usepackage{mathabx}
\usepackage[all,cmtip]{xy}
\usepackage{graphicx}
\usepackage{tikz-cd}
\usetikzlibrary{arrows.meta, backgrounds, positioning, calc}
\usepackage{color}
\usepackage{hyperref}
\usepackage{pgffor}
\usepackage{bbm}
\usepackage{siunitx}

\newtheorem{theorem}{Theorem}

\newtheorem{corollary}[theorem]{Corollary}
\newtheorem*{corollary*}{Corollary}

\theoremstyle{definition}

\newtheorem*{question}{Question}

\theoremstyle{remark}

\renewcommand{\phi}{\varphi}
\renewcommand{\theta}{\vartheta}
\renewcommand{\epsilon}{\varepsilon}

\newcommand{\NN}{\mathbb{N}}

\newcommand{\RR}{\mathbb{R}}

\newcommand{\id}{\operatorname{id}}

\title{Fundamental Groups of Permutation Complexes}
\author[D. Govc]{Dejan Govc}
\thanks{Supported by the Slovenian Research Agency program P1-0292.}
\address{University of Ljubljana, Faculty of Mathematics and Physics, Jadranska ulica 21, 1000 Ljubljana, Slovenija}
\email{dejan.govc@fmf.uni-lj.si, ORCID: \href{https://orcid.org/0000-0002-7156-7966}{0000-0002-7156-7966}}

\begin{document}

\begin{abstract}
	We present an elementary argument showing that fundamental groups of permutation complexes are free. Hence, $\RR P^2$ cannot be realized by a permutation complex, answering a question of Chach\'olski, Levi and Meshulam.
\end{abstract}

\maketitle

The complex of injective words $\operatorname{Inj}([n])$ is the semisimplicial complex whose facets are the $n!$ permutations in $S_n$. As shown by Farmer \cite{Farmer} and Bj\"orner--Wachs \cite{Bjorner-Wachs}, $\operatorname{Inj}([n])\simeq\bigvee_{d_n}S^{n-1}$, where $d_n$ is the number of derangements in $S_n$. In \cite{CLM}, Chach\'olski, Levi and Meshulam undertook the study of the more general notion of permutation complexes. Given a set of permutations $S\subseteq S_n$, the \emph{permutation complex} $X(S)$ is defined as the semisimplicial complex\footnote{In \cite[Definition 1.4]{CLM}, $X(S)$ is defined in terms of order complexes instead. This means that $X(S)$, as defined in \cite{CLM}, is the barycentric subdivision of $X(S)$ as defined in the present paper.} whose facets are given by the elements of $S$. Hence, a $k$-simplex in $X(S)$ is precisely a subword $s_1s_2\ldots s_{k+1}\leq\sigma\in S$. (For example, $X(1234,4312)$ has both $124$ and $412$ as $2$-simplices.) In particular, $X(S)$ is a directed simplicial complex in the sense of \cite[Section 4.1.2]{Frontiers}.

As shown in \cite{CLM}, permutation complexes have a rich homotopy theory. Specifically, by \cite[Theorem 1.5]{CLM}, every finite homotopy type $K$ can be realized stably by a permutation complex, i.e. there is an $r\in\NN$ and a permutation complex $X(S)$ such that $X(S)\simeq\Sigma^rK$. Therefore, it is natural to ask how many suspensions are necessary in such a representation. The present author found an example of $S\subseteq S_4$ consisting of $15$ permutations such that $X(S)\simeq\Sigma\RR P^2$ (see \cite[Paragraph following Example 2.6]{CLM}). In \cite[Concluding Remarks]{CLM} the following question was posed:
\begin{question}
Is there a set $S\subseteq S_n$ for some $n\in\NN$ such that $X(S)\simeq\RR P^2$?
\end{question}
We give a negative answer by showing that the fundamental group of a permutation complex is always free. The argument is elementary and uses the standard method to calculate the fundamental group from the CW structure (Hatcher \cite[Section 1.2]{Hatcher}). Let $n$ and $m$ be positive integers and let $X=X(S)$ be the permutation complex generated by $S=\{\sigma_1,\ldots,\sigma_m\}\subseteq S_n$. The main result is the following.

\begin{theorem}\label{thm:main}
	The fundamental group $G=\pi_1(X)$ is free.
\end{theorem}

\begin{corollary}
	There is no set $S\subseteq S_n$ for any $n\in\NN$ such that $X(S)\simeq\RR P^2$.
\end{corollary}

\enlargethispage{2\baselineskip}
\section{A Presentation of the Fundamental Group}

By renaming, if needed, we may assume that $\sigma_1=\id$. The $2$-skeleton $X(S)^{(2)}$ has a spanning tree $\Gamma=\{12,13,\ldots,1n\}$. Then $W=X(S)^{(2)}/\Gamma$ has:
\begin{itemize}
	\item A single $0$-cell.
	\item A $1$-cell $uv$ for every word $uv\leq\sigma\in S$ such that $u\neq 1$.
	\item A $2$-cell $uvw$ for every word $uvw\leq\sigma\in S$.
\end{itemize}
The presentation has generators $g_{uv}$ ($g_{1a}=1$) and relators $r_{uvw}=g_{uv}g_{vw}g_{uw}^{-1}$:
\[
\pi_1(X)=\langle E\mid R\rangle,\quad E=\{g_{uv}\mid uv\leq\sigma\in S, u\neq 1\},\quad R=\{r_{uvw}\mid uvw\leq\sigma\in S\}.
\]

\section{Proof of Theorem \ref{thm:main}}

Now let $x_a:=g_{a1}$. The relations involving vertex $1$ are as follows (if they occur):
\[
r_{1uv}=g_{uv},\qquad r_{u1v}=g_{u1}g_{uv}^{-1},\qquad r_{uv1}=g_{uv}g_{v1}g_{u1}^{-1}.
\]
For each $g_{uv}$ that occurs with $1\notin\{u,v\}$, there is at least one triangle of the form $1uv$, $u1v$ or $uv1$ that occurs in $X(S)$. Choose the first one that occurs (in the listed order), call it $t_{uv}$ and use the corresponding Tietze move to eliminate $g_{uv}$; i.e. remove the generator $g_{uv}$ and replace it by $\overline{g}_{uv}$ in each relation, where if $t_{uv}=1uv$, then $\overline{g}_{uv}=1$; if $t_{uv}=u1v$, then $\overline{g}_{uv}=x_u$; and if $t_{uv}=uv1$, then $\overline{g}_{uv}=x_u x_v^{-1}$.

The generators with $1\in\{u,v\}$ remain unchanged, i.e. $\overline{g}_{1a}=1$ and $\overline{g}_{a1}=x_a$.

This yields a new presentation with generators $x_a$ and relators $\overline{r}_{uvw}=\overline{g}_{uv}\overline{g}_{vw}\overline{g}_{uw}^{-1}$:
\[
\pi_1(X)=\langle T\mid \overline{R}\rangle,\quad T=\{x_a\mid a1\leq\sigma\in S, a\neq 1\},\quad \overline{R}=\{\overline{r}_{uvw}\mid uvw\leq\sigma\in S\}.
\]
We claim that we can eliminate every relation corresponding to a triangle not containing $1$. To do this, observe that any such triangle lies in a tetrahedron containing $1$, i.e. either $1uvw$, $u1vw$, $uv1w$ or $uvw1$. In each of these four cases, the relator $\overline{r}_{uvw}$ decomposes as a product of relators (or inverses/conjugates of relators) corresponding to the remaining three triangles. We treat each case separately:

{\bf Case 1:} $1uvw$. Note that $uv, vw$ and $uw$ lie in triangles $1uv,1vw$ and $1uw$ respectively, so the reduced relation is $\overline{r}_{uvw}=\overline{g}_{uv}\overline{g}_{vw}\overline{g}_{uw}^{-1}=1$, i.e. trivial.

{\bf Case 2:} $u1vw$. Since $uv, vw$ and $uw$ lie in $u1v,1vw$ and $u1w$ respectively, we have the relations $\overline{r}_{u1v}=x_u\overline{g}_{uv}^{-1}$, $\overline{r}_{1vw}=\overline{g}_{vw}=1$ and $\overline{r}_{u1w}=x_u\overline{g}_{uw}^{-1}$. The desired factorization is: $\overline{r}_{uvw}=\overline{g}_{uv}\overline{g}_{vw}\overline{g}_{uw}^{-1}=\overline{r}_{u1v}^{-1}\overline{r}_{u1w}$.

{\bf Case 3:} $uv1w$. Since $uv, vw$ and $uw$ lie in triangles $uv1,v1w$ and $u1w$ respectively, we have the relations $\overline{r}_{uv1}=\overline{g}_{uv}x_vx_u^{-1}$, $\overline{r}_{v1w}=x_v\overline{g}_{vw}^{-1}$ and $\overline{r}_{u1w}=x_u\overline{g}_{uw}^{-1}$. The desired factorization is $\overline{r}_{uvw}=\overline{g}_{uv}\overline{g}_{vw}\overline{g}_{uw}^{-1}=\overline{r}_{uv1}(x_ux_v^{-1}\overline{r}_{v1w}^{-1}x_vx_u^{-1})\overline{r}_{u1w}$.

{\bf Case 4:} $uvw1$. Since $uv, vw$ and $uw$ lie in triangles $uv1,vw1$ and $uw1$ respectively, we have the relations $\overline{r}_{uv1}=\overline{g}_{uv}x_vx_u^{-1}$, $\overline{r}_{vw1}=\overline{g}_{vw}x_wx_v^{-1}$ and $\overline{r}_{uw1}=\overline{g}_{uw}x_wx_u^{-1}$. Therefore, $\overline{r}_{uvw}=\overline{g}_{uv}\overline{g}_{vw}\overline{g}_{uw}^{-1}=\overline{r}_{uv1}(x_ux_v^{-1}\overline{r}_{vw1}x_vx_u^{-1})\overline{r}_{uw1}^{-1}$.

Hence, the presentation simplifies to:
\[
\pi_1(X)=\langle T\mid R'\rangle,\quad T=\{x_a\mid a1\leq\sigma\in S, a\neq 1\},\quad R'=\{\overline{r}_{uvw}\mid 1\leq uvw\leq\sigma\in S\}.
\]
Now observe that every relator $\overline{r}_{uvw}$ with $1\leq uvw$ is either trivial or reduces to $x_a=1$ or $x_a=x_b$. In detail: $\overline{r}_{1uv}=\overline{g}_{uv}=1$; $\overline{r}_{u1v}=x_u\overline{g}_{uv}^{-1}$, which is either $x_u$ or $1$, depending on $t_{uv}$; and $\overline{r}_{uv1}=\overline{g}_{uv}x_vx_u^{-1}$, which is either $x_vx_u^{-1}$, $x_ux_vx_u^{-1}$ or $1$, again depending on $t_{uv}$. Hence, this is the presentation of a free group.

\subsection*{Acknowledgments}

The author would like to thank Ran Levi and Assaf Libman for reading a preliminary version of the paper.

\subsection*{AI disclosure}

The proof of the main theorem was discovered using ChatGPT-6 Astra Max. An earlier proof\footnote{Presented by the author at Homology and Homotopy in the World of Graphs in Lausanne, Switzerland, on September 14, 2026 \cite{Govc-Lausanne}.} was discovered using ChatGPT-5.6 Sol Max, but is not included in the paper, since a simpler argument now exists. The author has independently verified the argument, wrote the final exposition and takes full responsibility for the correctness of this paper.

\bibliographystyle{amsplain}
\bibliography{paper}

\end{document}